\DeclareFontFamily{U}{wncy}{}
\DeclareFontShape{U}{wncy}{m}{n}{%
   <5>wncyr5%
   <6>wncyr6%
   <7>wnyr7%
   <8>wncyr8%
   <9>wncyr9%
   <10>wncyr10%
   <11>wncyr10%
   <12>wncyr6%
   <14>wncyr7%
   <17>wncyr8%
   <20>wncyr10%
   <25>wncyr10}{}
\DeclareMathAlphabet{\cyrille}{U}{wncy}{m}{n}
\documentclass[a4paper,11pt]{article}

\usepackage{latexsym}
\usepackage{amsmath,amsfonts,amssymb,amsthm}
\newtheorem{thm}{Theorem}[section]

\newtheorem{cor}[thm]{Corollary}

\newtheorem{rem}[thm]{Remark}

\input epsf

\title{On the enumeration of labelled hypertrees and of 
labelled bipartite trees}
\author{Roland Bacher}

\begin{document}
\maketitle
%\par polyuniv1.tex dans recherche/polyuniv 

{\sl Abstract\footnote{
Math. class: Primary: 05C30, Secondary: 05A05, 05A15, 05A19,
05C05, 05C07, 05C65}: We give a simple formula for the number of hypertrees
with $k$ hyperedges of given sizes and $n+1$ labelled vertices
with prescribed degrees. A slight generalization of this formula
counts labelled bipartite trees with prescribed degrees
in each class of vertices.
}

%fichier enumhyptree1.tex dans polyuniv

%\section{Introduction} 

%We describe a formula enumerating the number of hypertrees having labelled
%vertices with prescribed degrees and consisting of a given number
%of hyperedges with prescribed sizes. Our proof is bijective and
%is in some sense dual to the so-called Pr\"ufer code.

%Since hypergraphs are almost in one-to-one correspondence with
%bipartite (ordinary) graphs, a sligth variation of the above
%formula gives the number of labelled bipartite trees having 
%vertices with prescribed degrees.

\section{Main results}

\subsection{Labelled hypergraphs}

A \emph{(finite) hypergraph} is a pair $(\mathcal V,\mathcal E)$
consisting of a finite set $\mathcal V$ of \emph{vertices} 
and of a set $\mathcal E$ of \emph{hyperedges}
given by subsets of $\mathcal V$ containing at least two elements.
We define the \emph{size} of a hyperedge $E$ as 
the number $\mathrm{size}(E)$ of 
vertices contained in $E$ and the \emph{degree} $\deg v$
of a vertex
$v$ as the number of hyperedges containing $v$.
A vertex of degree $1$ is also called a \emph{leaf}.
The obvious linear relation
\begin{eqnarray}\label{linrelve}
\sum_{v\in\mathcal V}\deg(v)=\sum_{E\in\mathcal E}\mathrm{size}(E)
\end{eqnarray}
links the total sum of vertex-degrees to the total sum
of hyperedge-sizes.

Two distinct vertices $v,w\in\mathcal V$ are \emph{adjacent} 
or \emph{neighbours} if they are both contained in some  
hyperedge $E$ of $\mathcal E$.
A \emph{path} of \emph{length $k$} joining two vertices 
$v,w$ in a hypergraph $(\mathcal V,\mathcal E)$ is a sequence
$v_0=v,v_1,\dots,v_k=w$ involving only
consecutively adjacent vertices. 
A hypergraph is \emph{connected} if any pair of vertices 
can be joined by a path. 
The \emph{(combinatorial) distance} between two vertices $v,w$
of a connected hypergraph 
is the minimal length of a path in the set of all paths
joining $v$ and $w$.
A \emph{cycle} of length $k$ (for $k\geq 3$) is a closed 
path consisting of $k$ distinct vertices. 
A \emph{hyperforest} is a hypergraph $F$ such that two distinct 
hyperedges of $F$ intersect in at most a common vertex and such that
every cycle of $F$ is contained
in a hyperedge. A \emph{hypertree} is a connected hyperforest.
Induction on the number $k$ of hyperedges in a hypertree consisting of 
$n+1$ vertices shows the relation
\begin{eqnarray}\label{eqntreesize}
\sum_{E\in\mathcal E}\mathrm{size}(E)=n+k\ .
\end{eqnarray}

Let $\lambda$ be a partition of $n=\sum_{j=1}^k\lambda_j$ 
having exactly $k\leq n$
nonzero parts $\lambda_1\geq\lambda_2\geq\dots\geq \lambda_k\geq 1$.
%We denote by $\nu_1,\dots,\nu_n$ the multiplicities
%$\nu_j=\sharp\{i\ \vert\lambda_i=j\}$ counting parts of length $j$ 
%in $\lambda$.
Let $\mu=(\mu_0,\dots,\mu_n)\in\mathbb N^{n+1}$ 
be a vector with coefficients formed by $n+1$ natural integers 
$\mu_0,\dots,\mu_n$ summing up to $k-1=\sum_{j=0}^n\mu_j$.

\begin{thm}\label{thmenumtsizedeg} 
The number of hypertrees having $n+1$ vertices $\{0,\dots,
n\}$ of degrees $\deg(i)=1+\mu_i$ and $k$ hyperedges of sizes
$1+\lambda_1,\dots,1+\lambda_k$ is given by 
\begin{eqnarray}\label{formenumsizedeg}
{n\choose\lambda}\frac{1}{\prod_{j=1}^n(\nu_j)!}{k-1\choose
\mu}=
\frac{n!}{\left(\prod_{j=1}^k\lambda_j!\right)}
\frac{1}{\left(\prod_{j=1}^n(\nu_j)!\right)}
\frac{(k-1)!}{\left(\prod_{j=0}^n\mu_j!\right)}
\end{eqnarray}
%with ${n\choose \lambda}$ and ${k-1\choose \mu}$ denoting the obvious
%multinomial coefficients and 
with $\nu_j=\sharp\{i\ \vert\lambda_i=j\}$ 
counting parts of length $j$ in $\lambda$.
\end{thm}

Theorem \ref{thmenumtsizedeg} has the following equivalent formulation:

\begin{thm} \label{thmenumsize}
Denoting by $\mathcal{HT}_{\lambda}(n+1)$
the set of hypertrees with $n+1$ labelled vertices $\{0,\dots,n\}$
and with $k$ edges of size $1+\lambda_1,1+\lambda_2,\dots,1+\lambda_k$, 
we have
\begin{eqnarray}\label{formenumsize}
\sum_{T\in\mathcal{HT}_\lambda(n+1)}\prod_{j=0}^n
x_j^{\deg(j)}={n\choose\lambda}\frac{1}{\prod_{j=1}^n(\nu_j)!}
\left(x_0+x_1+\dots+x_n\right)^{k-1}\ .
\end{eqnarray}
\end{thm}

Equivalence between Theorem \ref{thmenumtsizedeg} and \ref{thmenumsize}
is given by the binomial Theorem.

The identities (\ref{linrelve}) and (\ref{eqntreesize}) show
that the conditions $\sum_{i=1}^k\lambda_k=n$ and 
$\sum_{i=0}^n\mu_i=k-1$ are necessary for the existence of
hypertrees with edge-sizes $1+\lambda_1,\dots,1+\lambda_k$ and
vertex-degrees $1+\mu_0,\dots,1+\mu_n$. Theorem 
\ref{thmenumtsizedeg} or \ref{thmenumsize} shows that 
they are also sufficient.

Multiplication by $(n+1)$ of the results given by Theorem
\ref{thmenumtsizedeg} and \ref{thmenumsize} gives enumerative
results for rooted labelled hypertrees.

Removal of the vertex $0$ in trees enumerated by Theorem 
\ref{thmenumtsizedeg} gives enumerative results for planted forests
with $\mu_0+1$ connected components (inducing an error of $1$ in
the degree of root vertices and in the size of hyperedges containing a root).

Counting trees or hypertrees with labelled vertices is a fairly old
sport and started with Sylvester \cite{Sy} and Cayley \cite{Ca}
(according to the notes of
Chapter 5 in \cite{St}) mentionning the total number 
$$(n+1)^{n-1}$$
of labelled trees on $n+1$ vertices. This corresponds of course 
to the specialization $x_0=\dots=x_n=1$ in Theorem \ref{thmenumsize}
for $\lambda=(1,1,\dots,1)$ the trivial partition of $n$.
In \cite{EE}, Erd\'ely and Etherington refined Cayley's theorem 
by enumerating labelled (ordinary)
trees with given vertex-degrees 
(corresponding to the case of the trivial partition 
$\lambda=(1,1,\dots,1)$ of $n$ in Theorem \ref{thmenumtsizedeg} or
\ref{thmenumsize}),
see also Theorem 5.3.10 in \cite{St}.
For trees, one can also consult the monograph \cite{Mo} or the
numerous more recent literature.

Denoting by $S_2(n,k)$ the Stirling number of the second kind enumerating 
partitions of $\{1,\dots,n\}$ into $k$ non-empty subsets, 
summing identity (\ref{formenumsize})
over all partitions of $\{1,\dots,n\}$ into exactly $k$ parts 
and setting $x_0=x_1=\dots=x_n=1$ yields the number
\begin{eqnarray}\label{thmHusimi}
&(n+1)^{k-1}S_2(n,k)
\end{eqnarray} 
of hypertrees with $n+1$ labelled vertices and
$k$ hyperedges given by Husimi in \cite{Hu}, see also 
\cite{Kr}, \cite{Ta}, \cite{GK} and \cite{BCS}
for other treatments and related results.

%Theorem \ref{thmenumtsizedeg} seems to have been overlooked although 
%there are some partial results in its direction, reformulated 
%in the language of planted labelled forests.

It is perhaps worthwile to mention the following two corollaries of 
Theorem \ref{thmenumtsizedeg}:

The first result counts weighted hypertrees and is
a kind of counterpart of Husimi's result (\ref{thmHusimi})
in the sense that it involves Stirling numbers
of the first kind counting the number $(-1)^{n+k}S_1(n,k)$ of permutations 
of $\{1,\dots,n\}$ involving $k$ disjoint cycles:

\begin{cor} We have 
$$\sum_{T\in\mathcal{HT}_k(n+1)}w(T)=(-1)^{n+k}(n+1)^{k-1}S_1(n,k)$$
or more precisely
$$\sum_{T\in\mathcal{HT}_k(n+1)}w(T)
\prod_{j=0}^n x_j^{\deg(j)}=(-1)^{n+k}(x_0+x_1+\dots+x_n)^{k-1}S_1(n,k)$$
where $\mathcal{HT}_k(n+1)$ denotes the set of all labelled hypertrees
with $k$ hyperedges and vertices $\{0,\dots,n\}$, where $w(T)=\prod_{j=1}^k
(\lambda_j-1)!$ for a labelled hypertree $T$ with $k$ hyperedges
of size $1+\lambda_1,\dots,1+\lambda_k$ and where $S_1(n,k)$
is the Stirling number of the first kind defined by 
$\sum_{k=0}^nS_1(n,k)x^k=\prod_{j=0}^{n-1}(x-j)$.
\end{cor}

{\bf Proof} Observe that 
$${n\choose\lambda}\frac{1}{\prod_{j=1}^n(\nu_j)!}\prod_{j=1}^k(\lambda_j-1)!$$
counts the number of permutations of $\{1,\dots,n\}$
in the conjugacy class consisting of products of $k$ disjoint 
cycles with lengths $\lambda_1,\dots,\lambda_k$.
Summing over all partitions of $\{1,\dots,n\}$ into $k$ non-empty subsets
and applying Theorem \ref{thmenumsize} (with $x_0=\dots=x_n=1$ 
for the first formula)
yields the result since $(-1)^{n+k}S_1(n,k)$ counts
the total number of permutations of $\{1,\dots,n\}$
consisting of $k$ disjoint cycles.\hfill
$\Box$

\begin{cor} Let $\mathcal{HT}_k(n+1)$ be the set of all labelled 
hypertrees with $k$ hyperedges and $n+1$ vertices.
The two random variables given by the sizes of hyperedges and 
by the degrees of vertices are independent for a random hypertree $T$
choosen with uniform probability in $\mathcal{HT}_k(n+1)$.

More precisely, a random hypertree, choosen with uniform probability 
among all labelled hypertrees with $k$ hyperedges and vertices $\{0,
\dots,n\}$, has edge-sizes $1+\lambda_1,\dots,1+\lambda_k$
associated to a partition $\lambda$ of $n$ with probability
$$\frac{1}{S_2(n,k)}{n\choose \lambda}\frac{1}{\prod_{j=1}^k
(\nu_j)!}$$
with $S_2(n,k)$ denoting the Stirling number of the second kind 
counting partitions of $\{1,\dots,n\}$ into $k$ non-empty subsets
and with $\nu_j=\sharp\{i\ \vert \ \lambda_i=j\}$ counting
the number of parts equal to $j$ in $\lambda$.

Such a random tree has vertices $0,\dots,n$ 
of degrees $1+\mu_0,\dots,1+\mu_n$
with probability $\frac{1}{(n+1)^{k-1}}{k-1\choose \mu}$.
\end{cor}

{\bf Proof} This is an immediate consequence of the fact that the right 
side of formula~(\ref{formenumsizedeg}) factors into a product of 
two terms depending only on hyperedge-sizes, respectively
vertex-degrees.\hfill$\Box$

%%%%%%%%%%%%%%%%%%%%%%%%%%%%%%%%

\subsection{Labelled bipartite trees}

A bipartite graph is an ordinary graph with vertices 
$\mathcal V=\mathcal V_1\cup
\mathcal V_2$ partitioned into two subsets such that no pair of adjacent
vertices is in the same class.

Hypergraphs are in one-to-one correspondence with certain bipartite
graphs as follows (see for example page 5 of \cite{Bo}):
To a hypergraph $(\mathcal V,\mathcal E)$ we
associate the bipartite graph with vertices in the first
class representing elements of $\mathcal V$, vertices of the second
class representing elements of $\mathcal E$ and with edges
encoding incidence (ie. there is an egde relating a vertex $v$ to
a hyperedge $E$ if $v$ belongs to $E$).
Vertices of edge-type (representing elements of $\mathcal E$) 
have to be of degree at least $2$ and every bipartite graph
having only vertices of degree at least $2$ in its second class
of vertices corresponds to a hypergraph. Hypertrees are encoded by (ordinary) 
trees with no leaves in their second bipartite class of vertices.

\begin{thm}\label{thmbipartite} 
Given two natural integers $a,b$ and two 
integral vectors $\alpha=(\alpha_0,\dots,\alpha_a)
\in\mathbb N^{a+1}$ and $\beta=
(\beta_0,\dots,\beta_b)$ such that $b=\sum_{i=0}^a\alpha_i$
and $a=\sum_{i=0}^b\beta_i$, the number of labelled bipartite 
trees having vertex bipartition $\mathcal U\cup \mathcal V$
with vertices $\mathcal U=\{u_0,\dots,u_a\}$ of degree 
$\deg u_i=1+\alpha_i$ and vertices 
$\mathcal V=(v_0,\dots,v_b)$ of degree $\deg v_i=1+\beta_i$
is given by
\begin{eqnarray}\label{formulebipartite}
{a\choose \beta}{b\choose \alpha}=
\frac{a!\ b!}{\left(\prod_{i=0}^b\beta_i!\right)
\left(\prod_{i=0}^a\alpha_i!\right)}\ .
\end{eqnarray}
Equivalently, we have
\begin{eqnarray}\label{formulebipxy}
\sum_{T\in\mathcal T(a+1,b+1)}\left(\prod_{i=0}^a
x_i^{\deg(v_i)-1}\right)\left(\prod_{i=0}^b
y_i^{\deg(u_i)-1}\right)=\left(x_0+\dots+x_a\right)^b\left(
y_0+\dots+y_b\right)^a
\end{eqnarray}
where $\mathcal T(a+1,b+1)$ denotes the set of labelled
trees with $a+1$ vertices in the first class and $b+1$ vertices
in the second class of its vertex-partition.
%In particular, the set $\mathcal T(a+1,b+1)$ consists of 
%$(a+1)^b(b+1)^a$ elements.
\end{thm}

Setting $x_0=\dots=x_a=y_0=\dots=y_b=1$ in Formula 
(\ref{formulebipxy}) yields the well-known number
$(a+1)^b(b+1)^a$ of spanning trees in the complete 
bipartite graph $K_{a+1,b+1}$.

\begin{rem} Theorem \ref{thmbipartite} implies that a 
bipartite random tree (choosen with uniform probability)
with vertex-bipartition $\{u_0,\dots,u_a\}\cup\{v_0,\dots,v_b\}$
has vertices of degree $\deg(u_i)=1+\alpha_i$ for $i=0,\dots,a$
with probability $\frac{1}{(1+a)^b}{b\choose \alpha}$ independently
of the degrees of the vertices $v_0,\dots,v_b$.
\end{rem}

We give two proofs of Theorem \ref{thmbipartite}. The first proof shows
that it is essentially equivalent to Theorem \ref{thmenumtsizedeg}.
The second proof is bijective: We construct a map
$T\longmapsto (W(T),W'(T))$ from trees with vertex-bipartition 
$\mathcal U\times \mathcal V$ into $\mathcal V^{\sharp(\mathcal U)-1}
\times \mathcal U^{\sharp(\mathcal V)-1}$ which is one-to-one and
respects degrees: 
A vertex $u\in\mathcal U$, respectively $v\in\mathcal V$, of $T$ 
is involved with multiplicity $\deg(u)-1$ in $W'(T)$, respectively
with multiplicity $\deg(v)-1$ in $W(T)$.

%%%%%%%%%%%%%%%%%%%%%%%%%%%%%%%%%%%%%%%%%%%%%%%%%%%%%%%%%%%%%%%%%%%%%%

\section{Proof of Theorem \ref{thmenumtsizedeg}}
We give a bijective proof of Theorem \ref{thmenumtsizedeg}.
More precisely, we construct a map $T\longmapsto (\mathcal P(T),
W(T))$ which associates to a hypertree with $k$ hyperedges
of size $1+\lambda_1,\dots,1+\lambda_k$ and vertices $\{0,\dots,n\}$
of degrees $\mu_0,\dots,\mu_n$
a pair $(\mathcal P(T),W(T))$ formed by a partition 
$\mathcal P(T)$ of $\{1,\dots,n\}$
into $k$ subsets of cardinalities $\lambda_1,\dots,\lambda_k$
and by a word $W(T)\in\{0,\dots,n\}^{k-1}$ of length $k-1$ 
involving a letter $i$
of the alphabet $\{0,\dots,n\}$ with multiplicity $\mu_i$.
Theorem~\ref{thmenumtsizedeg} follows from the fact that the map
$T\longmapsto(\mathcal P(T),W(T)$ is one-to-one, since there are 
${n\choose\lambda}\frac{1}{\prod_{j=1}^n(\nu_j)!}$
possibilities for $\mathcal P(T)$ and ${k-1\choose\mu}$
possibilities for $W(T)$.
%%%%%%%%%%%%%%%%%%%%%%%%%
\subsection{The map $T\longmapsto \mathcal P(T)$}

A hyperedge $E$ of a hypertree with vertices $\{0,\dots,n\}$
contains a unique vertex $m(E)$ at minimal distance to
the vertex $0$. In particular, we have $m(E)=0$ 
if $E$ contains $0$. We call $m(E)$ the \emph{marked vertex of $E$}. 

Removing the marked vertex $m(E)$ from every hyperedge $E\in\mathcal E(T)$
of a hypertree $T$ with vertices $\{0,\dots,n\}$ yields a partition 
$\mathcal P(T)$
of $\{1,\dots,n\}$ with parts $E\setminus\{m(E)\}$ of cardinalities
$\left(\mathrm{size}(E)-1\right)$ indexed by the set $\mathcal E(T)$ of 
all hyperedges in $T$.

%%%%%%%%%%%%%%%%%%%%%%%%%
\subsection{Construction of $T\longmapsto W(T)$}

Given a hypertree $T$ with 
vertices $\{0,\dots,n\}$ and $k$ hyperedges,
we construct recursively a word $W(T)\in \{0,\dots,n\}^{k-1}$ 
which encodes exactly
the loss of information induced by the map 
$T\longmapsto \mathcal P(T)$.

The construction of the word $W(T)$ is somehow dual to the Pr\"ufer
code encoding labelled planted forests (see for example the first proof
of Theorem 5.3.2 in \cite{St}). Pr\"ufer codes
are defined by keeping track of neighbours of successively removed  
largest leaves in ordinary trees, 
the construction of the word $W(T)$
is based on local simplifications 
around largest non-leaves in hypertrees.

We start with a few useful definitions and notations:

A \emph{hyperstar} is a hypertree 
containing a vertex $v$, called a \emph{center} of the hyperstar,
at distance at most $1$ from all other vertices.
A center of a hyperstar is unique (and given by the intersection of
two arbitrary hyperedges) except in the degenerate case
where the hyperstar consists of a unique hyperedge.

As above, we denote by $m(E)$ the marked vertex realizing the distance
to $0$ of a hyperedge $E$.
The remaining vertices of $E$ are \emph{unmarked} vertices.

Given $i$ in $\{1,\dots,n\}$, we denote by $U(i)\in \mathcal E(T)$ the unique
hyperedge of $T$ containing $i$ as an unmarked vertex. Similarly,
we denote by $P(i)$ the unique set $U(i)\setminus m(U(i))$ 
containing $i$ of the partition $\mathcal P(T)$. The set $P(i)$
can be constructed by removing the marked vertex from the unique 
hyperedge $U(i)$ containing $i$ as an unmarked vertex. 

We associate to every hypertree $T$ with $k$ hyperedges and
vertices $\{0,\dots,n\}$ of degrees $1+\mu_0,\dots,1+\mu_n$
a word $W(T)$ of $\{0,\dots,n\}^{k-1}$ involving
$\mu_i$ copies of a letter $i$ in 
$\{0,\dots,n\}$. The word $W(T)=w_1\dots w_{k-1}$ is defined 
recursively as follows:

We set $W(T)=0^{k-1}$ if $T$ is a hyperstar centered at $0$.

Otherwise, there exists a largest integer $a$ in $\{1,\dots,n\}$
such that $\mu_a>0$ and $\mu_{a+1}=\mu_{a+2}=\dots=\mu_n=0$.
We denote by $A_1,\dots,A_{k-1}\subset\{1,\dots,n\}\setminus\{a\}$ 
all $k-1$ elements
of $\mathcal P(T)\setminus \{P(a)\}$ not containing $a$ with 
indices determined by requiring $\min A_i<\min A_j$ if $i<j$.
Exactly $\mu_a$ elements among $A_1,\dots,A_{k-1}$
correspond to the set $m^{-1}(a)$ of hyperedges with marked vertex $a$.
Let $\alpha_1,\dots,\alpha_{\mu_a}$ be the associated indices.
Otherwise stated, $A_{\alpha_j}\cup\{a\}$
is a hyperedge with marked vertex $a$ of $T$ for $j=1,\dots,\mu_a$.
The indices $\alpha_1,\dots,\alpha_{\mu_a}\in\{1,\dots,k-1\}$
define the $\mu_a$ letters $w_{\alpha_1}=\dots=w_{\alpha_{\mu_a}}=a$ 
of the word $W(T)=w_1w_2\dots w_{k-1}$.
Removing these $\mu_a$ letters from the word $W(T)$ leaves a word
$W'$ which we define recursively by the identity 
$W'=W(T_a)$ 
where $T_a$ is obtained from $T$ by merging all hyperedges of $T$ containing
$a$ into a unique hyperedge. 
More precisely, $T_a$ is constructed by
removing first all hyperedges containing $a$ from of $T$,
followed by the adjunction of one new hyperedge 
consisting of $a$ and of all its neighbours in $T$.
The hypertree $T_a$ has strictly fewer hyperedges than $T$.
All vertices except $a$ have the same degree in $T$ and in $T_a$ and 
$a$ is a leaf in $T_a$.

Since all vertices $a,a+1,\dots,n$ are leaves of $T_a$,
we have $W'\in\{0,1,\dots,a-1\}^{k-1-\mu_a}$ and this processus
stops eventually.

Remark that the $\mu_a$ positions of all letters equal to $a$ in $W(T)$ 
encode exactly the elements of $\mathcal P(T)$ associated to hyperedges 
with marked vertex $a$ in a hypertree $T$. 
This implies that the map $T\longmapsto
(\mathcal P(T),W(T))$ is into.

%%%%%%%%%%%%%%%%%%%%%%%%%%%%%%

\subsection{Construction of the reciprocal map $(\mathcal P,W)\longmapsto 
T$}

We claim that the map $T\longmapsto (\mathcal P(T),W(T))$ is
one-to-one: Indeed, let $\mathcal P$ 
be a partition of $\{1,\dots,n\}$ into $k$ non-empty subsets and let
$W\in\{0,\dots,n\}^{k-1}$ be a word of length $k-1$ with letters
in the alphabet $\{0,\dots,n\}$. 
If $W=0^{k-1}$, the pair $(\mathcal P,W)$ corresponds to the
hyperstar centered at $0$ with hyperedges obtained by adding 
the central vertex $0$ to every subset $A\subset\{1,\dots,n\}$
involved in the partition $\mathcal P$.

Otherwise, let $a\in\{1,\dots,n\}$ be the largest strictly positive integer 
involved with strictly positive multiplicity $\mu_a>0$ 
in $W=w_1\dots w_{k-1}$ and let 
$1\leq \alpha_1<\dots<\alpha_{\mu_a}\leq k-1$ denote the $\mu_a$ indices 
defined by $w_{\alpha_1}=\dots=w_{\alpha_{\mu_a}}=a$.
We denote by $P(a)\in \mathcal P$ the unique subset of 
$\mathcal P$ containing the element $a$.
Let $A_1,\dots,A_{k-1}$ be the $k-1$ remaining elements of $\mathcal P$
corresponding to subsets of $\{1,\dots,\dots,n\}\setminus\{a\}$.
Indices of $A_1,\dots,A_{k-1}$
are defined by the requirement $\min A_i<\min A_j$ if $i<j$.
The sets $A_{\alpha_j}\cup\{a\}$, $j=1,\dots,\mu_a$, 
are then by construction the $\mu_a$ hyperedges 
with marked vertex $a$ of a
tree $T$ such that $\mathcal P=\mathcal P(T)$ and 
$W=W(T)$. The remaining hyperedges of such a tree $T$ are defined as follows:
merge the $\mu_a+1$ elements $P(a),A_{\alpha_1},\dots,A_{\alpha_{\mu_a}}$
of $\mathcal P$
into a unique subset $A=P(a)\cup\bigcup_{j=1}^{\mu_a}A_{\alpha_j}$ 
of $\{1,\dots,n\}$ and complete
this subset to a partition $\mathcal P'$ of $\{1,\dots,n\}$
by adjoining all elements of $\mathcal P$ not contained in $A$.
Similarly, define a word $W'$ obtained from $W$ by removing all $\mu_a$ 
occurences of the letter $a$. 
The pair $(\mathcal P',W')$ defines then recursively
a hypertree $T'$. Hyperedges of the tree $T$ not containing $a$ are then
given by hyperedges of $T'$ not containing $a$. The unique hyperedge $U(a)$
of $T$ containing the vertex $a$ as an unmarked vertex is obtained 
by adjoining to the set $P(a)$ the marked vertex of the unique hyperedge
in $T'$ with unmarked vertices given by $A=P(a)\cup\bigcup_{j=1}^{\mu_a}
A_{\alpha_j}\in\mathcal P'$.
Remark that the tree $T'$ is simpler than the final tree $T$ in the sense
that $a$ is a non-leaf of $T$ (the vertices $a+1,a+2,\dots,n$ are however
leaves of $T$), but is a leaf (together with $a+1,\dots,n$) of the tree $T'$.
Thus the construction stops eventually. 

A tree $T$ constructed in this way is the unique tree satisfying
$\mathcal P=\mathcal P(T)$ and $W=W(T)$. 
The map $T\longmapsto (\mathcal P(T),W(T))$ is thus also onto.
This ends the proof of Theorem \ref{thmenumtsizedeg}. \hfill$\Box$

\begin{rem} The bijection $T\longmapsto (\mathcal P(T),W(T))$
is not completely natural 
in the sense that it depends on the choice of a particular 
vertex $v$ (given by the vertex $0$ in our case), on the choice
of a linear order of the remaining vertices and on the
choice of a suitable order relation for subsets of partitions of 
$\mathcal V\setminus \{v\}$.  
\end{rem}

\begin{rem} The action of the symmetric group $\mathcal S_n$
on partitions of $\{1,\dots,n\}$ and the action of $\mathcal S_{k-1}$ 
permuting letters in a word of length $k-1$ induce a transitive action
of $\mathcal S_n\times \mathcal S_{k-1}$
on labelled
hypertrees with $k$ hyperedges of given sizes 
and vertices $0,\dots,n$ of given degrees.
\end{rem}

%%%%%%%%%%%%%%%%%%%%%%%%%%%%%%%%%%%%%%%%%%%%%%%%%%%%%%%%%%%%%%%%%%%

\section{Proofs of Theorem \ref{thmbipartite}}

We give two proofs of Theorem \ref{thmbipartite}.
The first proof consists in showing that it is essentially
equivalent to Theorem \ref{thmenumtsizedeg}.
The second proof is obtained by a minor modification of the bijective proof 
given above for Theorem \ref{thmenumtsizedeg}.

{\bf First proof} Suppose first that $\beta_0,\dots,\beta_b$ are all
strictly positive. We consider thus bipartite graphs having a
vertex-bipartition
$\mathcal U\cup\mathcal V$ with no leaves in $\mathcal V$.
Interpreting vertices of $\mathcal V$ as hyperedges, ordered
by size, the number of such graphs is obtained by multiplying the corresponding
number of labelled hypertree (given by formula (\ref{formenumsizedeg})
with $n=a,\lambda=\{\beta_0,\dots,\beta_b\}$ in decreasing order,
$k=b+1$ and $\mu=\{\alpha_0,\dots,\alpha_a\}$ in decreasing order)
by $\prod_{j=1}^b\nu_j!$
where $\nu_j=\sharp\{i\ \vert\ \beta_i=j\}$ counts the number of 
vertices of degree $j+1$ in $\mathcal V$.
We get thus in this case the equivalence between 
Theorem \ref{thmbipartite} and Theorem \ref{thmenumtsizedeg}.
The general case is by induction on the number of leaves in $\mathcal V$.
Indeed, the last such leave can be adjacent to any non-leaf in $\mathcal U$
and we get thus the recursion
$$\sum_{k,\alpha_k\geq 1}{a\choose \beta}{b-1\choose 
\alpha_0,\dots,\alpha_{k-1},\alpha_k-1,\alpha_{k+1},\dots,\alpha_a}=
{a\choose \beta}{b\choose \alpha}$$
for the total number of possible bipartite trees.\hfill$\Box$

{\bf Second proof}
We construct a map which associates to a tree $T$ having bipartite vertices
$\mathcal U\cup \mathcal V$, with vertices
$\mathcal U=\{u_0,\dots,u_\alpha\}$ of degrees
$\alpha_0,\dots,\alpha_a$ in  the first class 
and vertices $\mathcal V=\{v_0,\dots,v_\beta\}$
of degrees $\beta_0,\dots,\beta_b$ in the second class, a word 
$(W,W')=(W(T),W'(T))\in \mathcal V^a\times \mathcal U^b$
such that a vertex $v_i$ is involved $\beta_i$ times in $W$ and 
a vertex $u_i$ is involved $\alpha_i$
times in $W'$.
This implies the result since there are ${a\choose \beta}$ possibilities
for $W$ and ${b\choose \alpha}$ possibilities for $W'$.

We root a tree $T$ with vertex-bipartition 
$\mathcal U\times \mathcal V$ as above at the vertex $u_0$ of $\mathcal U$
and we orient all edges of $T$ away from the root vertex $u_0$. 
An edge joining two neighbouring vertices $s,t$ with $s$ 
closer to $u_0$ than $t$ is thus oriented from $s$ to $t$. We 
call the vertex $s$ of such an edge the \emph{parent} of $t$ and we write 
$s=p(t)$. Similarly, we call $t$ a \emph{child} of $s$.
Every vertex other than $u_0$ 
has a unique parent and edges of $T$ are in bijection with 
$\{u_1,\dots,u_a\}\cup \mathcal V$ by considering the unique edge 
joining a vertex $s\not=u_0$ to its parent $p(s)$.

We consider the total order induced by indices on 
both sets $\mathcal U$ and $\mathcal V$.

The word $W$ encodes all edges of $T$ starting at a vertex of $\mathcal V$
and ending at a vertex in $\{u_1,\dots,u_a\}$.
More precisely, $W$ is given by the word 
$$p(u_1)p(u_2)\dots p(u_a)$$
encoding the parents of $u_1,u_2,\dots,u_a$. Since every element $v_i$ in 
$\mathcal V$ is the parent of exactly $\beta_i$ vertices in $\mathcal U$,
the word $W$ involves a vertex $v_i\in\mathcal V$ with multiplicity $\beta_i$.

The recursive definition of the word $W'$ encoding all edgeds starting 
at an element of $\mathcal U$ and ending at an element of $\mathcal V$
is more involved. More precisely, for $i\geq 1$,
the positions of a letter $u_i$ in $W'$
encode the $\alpha_i$ edges of the form $\{v,u_i\}$ with $p(v)=u_i$ in 
the following way:

Suppose that for some integer $c\in\{1,\dots,a\}$ all positions 
of the letters $u_{c+1},u_{c+2},\dots,u_a$ in $W'$ 
are known. The $\alpha_{c+1}
+\alpha_{c+2}+\dots+\alpha_a$ positions of the letters $u_{c+1},
u_{c+2},\dots,u_a$ in $W'$ and the word $W$ define a subforest $F_c$ of $T$
consisting of all edges involving a vertex in $\{u_{c+1},u_{c+2},
\dots,u_a\}$ together with the $c$ edges of the form $\{u_i,p(u_i)\}$
joining the vertices $u_1,\dots,u_c$ to their predecessors 
$p(u_1),\dots,p(u_c)$ in $\mathcal V$. 
The forest $F_c$ has $u_0$ as an isolated 
vertex and contains 
$$b+1-(\alpha_{c+1}+\alpha_{c+2}+\dots+\alpha_a)
=1+\alpha_0+\alpha_1+\dots+\alpha_c$$
other connected components which are all rooted at a an element 
of $\mathcal V$ having no parent in its connected component of $F_c$.
We denote by $\tilde F_c$ the subforest of $F_c$ obtained by removing 
the isolated vertex $u_0$ and the connected component containing $u_c$ 
from $F_c$. The forest $\tilde F_c$ has exactly $\alpha_0+\dots+\alpha_c$
connected components which we order totally accordingly to the total order 
of the corresponding root-vertices.
The $\alpha_c$ children of $u_c$ are roots of $\alpha_c$
connected components of $\tilde F_c$. 
The relative positions of these $\alpha_c$
connected components among the totally ordered set of all 
connected components of $\tilde F_c$, 
determine the $\alpha_c$ positions of the letter $u_c$ among 
the $\alpha_0+\dots+\alpha_c$ letters of $W'$ which are different from
$u_{c+1},u_{c+2},\dots,u_a$.

If all letters except $u_0$ of $W'$ are known, we complete $W'$
with $\alpha_0$ copies of the letter $u_0$ in the unique possible way.

It is easy to see that the map $T\longmapsto (W(T),W'(T))$
is into. We leave it to the reader to check that it is onto
by working out the obvious reciprocal map $(W,W')\longmapsto T$
(only the last step involving the tree $F_0$ is not contained in the
above description: one gets $T$ from the forest $F_0$ by joining 
$u_0$ to the root-vertices of the remaining $\alpha_0+1$ connected
components of $F_0$).
It defines thus a one-to-one map finishing the proof of
Theorem \ref{thmbipartite}.\hfill$\Box$
 
\begin{rem} (i) The map $T\longmapsto (W(T),W'(T))$ depends on
total orders of the vertices. It depends also an order 
on the sets $\mathcal U,\mathcal V$ involved in the vertex-bipartition
$\mathcal U\cup\mathcal V$.
This map is of course very similar to the map considered
in the proof of Theorem \ref{thmenumtsizedeg}.
There is however a subtle difference: The forests $F_c$
used in the recursive construction of $W'(T)$ are naturally rooted
(by the choice of the root-vertex $u_0$ in $\mathcal U$). 
This is not the case by its counterpart given by
the partition $\mathcal P'$ in the proof for hypertrees.

(ii) The obvious action (by permuting the positions of the letters
in $W$ and $W'$) of $\mathcal S_a\times \mathcal S_b$ induces a 
transitive action on the set of labelled bipartite trees 
enumerated by formula (\ref{formulebipartite}).

(iii) Combining the Pr\"ufer code with the map $T\longmapsto
(W(T),W'(T))$ yields a one-to-one map
between $\{0,\dots,n\}^{n-1}$ and
$$\bigcup_{A
\subsetneq\{1,\dots,n\}}\{\{0\}\cup A\}^{n-1-\sharp(A)}\times
\{\{1,\dots,n\}\setminus A\}^{\sharp(A)}$$
(where we agree that the first vertex of a tree belongs always to 
the first subset of the vertex-bipartition)
illustrating the identity 
$$\left(x_0+\dots+x_n\right)^{n-1}=\sum_{A\subsetneq\{1,\dots,n\}}
\left(x_0+\sum_{j\in A}x_j\right)^{n-1-\sharp(A)}\left(
\sum_{j\in\{1,\dots,n\}\setminus A}x_j\right)^{\sharp(A)}$$
(for $n\geq 1$) having the specialization
$$(n+1)^{n-1}=\sum_{k=0}^{n-1}{n\choose k}(k+1)^{n-1-k}(n-k)^k\ ,$$
see \cite{Po} for similar identities.
Is there an easier way to describe such a bijection?
\end{rem}

\noindent Roland BACHER, Universit\'e Grenoble I, CNRS UMR 5582, Institut 
Fourier, 100 rue des maths, BP 74, F-38402 St. Martin d'H\`eres, France.

\noindent e-mail: Roland.Bacher@ujf-grenoble.fr


\begin{thebibliography}{99}

\bibitem{BCS} A. Bedini, S. Caracciolo, A. Sportiello, \lq\lq Hyperforests on
the Complete Hypergraph by Grassmann Integral Representation'', 
{\it J. Phys. A} {\bf 41} (2008), 205003, 28 pp.

\bibitem{Bo} B. Bollob\'as, {\it Graph theory. An introductory  course}, 
Graduate Texts in Mathematics, New York-Berlin, Springer, 1979.

\bibitem{Ca} A. Cayley, \lq\lq A theorem on trees'', {\it Quart. J. Math.} 
{\bf 23} (1889), 376--378; {\it Collected Papers}, vol. 13, Cambridge 
University Press, 1897, pp. 26--28.

\bibitem{EE} A. Erd\'ely, I.M.H. Etherington, \lq\lq Some problems of 
non-associative combinations (2)'', {\it Edinburgh Math. Notes}
{\bf 32} (1940), 7--12.

\bibitem{GK} I.M. Gessel, L.H. Kalikow, \lq\lq
Hypergraphs and a functional equation of Bouwkamp and de Bruijn'', 
{\it J. Combin. Theory Ser. A 110} (2005), no. 2, 275--289. 

\bibitem{Hu} K. Husimi, \lq\lq Note on Mayer's theory of cluster integrals'',
{\it Journal of Chemical Physics} {\bf 18} (1950), 682--684.

\bibitem{Kr} Kreweras, \lq\lq Counting problems in dendroids'' in 
{\it Combinatorial Structures and 
their Applications}, R.K. Guy et al. editors, 223-226, Gordon and Breach,
New York (1970).

\bibitem{Mo} J.W. Moon, {\it Counting labelled trees},
Canadian Mathematical Monographs, No. 1.

\bibitem{Po} T.D. Porter, \lq\lq Binomial Identities Generated by Counting
Spanning Trees'', {\it Ars Combinatoria} {\bf 82} (2007), 159--163.

\bibitem{St} R. P. Stanley, {\it Enumerative Combinatorics, Volume 2}, 
Cambridge University Press, Cambridge, 1999.

\bibitem{Sy} J.J. Sylvester, \lq\lq On the change of systems of 
independent variables'', {\it Quart. J. Math.} {\bf 1} (1857), 42--56; 
{\it Collected Mathematical Papers}, vol. 2, Cambridge, 1908, pp. 65--85.

\bibitem{Ta} L. Tak\`acs, \lq\lq On the number of distinct forests'',
{\it SIAM J.Disc. Math.} {\bf 3} (1990), 574--581.

\end{thebibliography}
\end{document}